\documentclass[12pt,reqno,a4paper]{amsart}
\usepackage[margin=1in]{geometry}
\usepackage{graphicx}
\usepackage{tikz}

\usepackage{amsmath,amsfonts,amssymb,amsthm,mathtools}
\usepackage[hidelinks,hyperindex,breaklinks]{hyperref}
\usepackage{graphics,xcolor,cleveref,esint,comment}
\usepackage{enumerate}

\newcommand{\R}{\mathbb{R}}

\DeclareMathOperator{\divr}{div}

\def\Xint#1{\mathchoice
   {\XXint\displaystyle\textstyle{#1}}%
   {\XXint\textstyle\scriptstyle{#1}}%
   {\XXint\scriptstyle\scriptscriptstyle{#1}}%
   {\XXint\scriptscriptstyle\scriptscriptstyle{#1}}%
   \!\int}
\def\XXint#1#2#3{{\setbox0=\hbox{$#1{#2#3}{\int}$}
     \vcenter{\hbox{$#2#3$}}\kern-.5\wd0}}

\def\dashint{\Xint-}

\newtheorem{Theorem}{Theorem}
\newtheorem{Remark}{Remark}

\newtheorem{Lemma}{Lemma}

\begin{document}
    \title{Scale invariant bounds for the Kelvin-Helmholtz instability}

    \author[K.~Kalinin]{Konstantin Kalinin}
    \address{Max Planck Institute for Mathematics in the Sciences, Inselstr. 22, 04103 Leipzig, Germany.}
    \email{konstantin.kalinin@mis.mpg.de}
    
    \author[G.~Menon]{Govind Menon}
    \address{Division of Applied Mathematics, Box F, Brown University, Providence, RI 02912, USA.}
    \email{govind\_menon@brown.edu}

    \author[B.~Wu]{Bian Wu}
    \address{Max Planck Institute for Mathematics in the Sciences, Inselstr. 22, 04103 Leipzig, Germany.}
    \email{bian.wu@mis.mpg.de}

    \address{Corresponding author: Konstantin Kalinin}

    \begin{abstract}
        We derive robust long-time a-priori estimates for the Navier-Stokes equation in a two-dimensional infinite strip which are uniform in the Reynolds number. These estimates provide several new scale invariant upper bounds for the size of the mixing layer in the Kelvin-Helmholtz instability.
    \end{abstract}
    \maketitle

    {\centering\footnotesize \em In memory of Charlie Doering.\par}

    \section{Introduction}
    %\subsection{Outline}

    \subsection{Overview}
    Our purpose in this paper is to provide scale-invariant upper bounds on the size of the mixing layer that arises due to the Kelvin-Helmholtz instability. This instability arises in various physical contexts, such as plasma physics \cite{johnson2014kelvin} and ocean mixing \cite{smyth2012ocean}. It is also of intrinsic mathematical interest.

    We model the Kelvin-Helmholtz instability by considering the development of a {\em viscous\/} incompressible flow, governed by the 2D Navier-Stokes equations, from an initial condition that is close to a discontinuous shear flow. The geometry of the flow is illustrated in Figure~\ref{fig:kh-instability} and the equations and initial conditions are described precisely in Section~\ref{subsec:equations} below. Our main result is a rigorous upper bound on the size of the mixing layer that is linear in time and independent of the Reynolds number, under the assumption that {\em the initial data has positive vorticity\/}. This result is motivated by both the rigorous analysis of the Euler equations, as well as recent work in the turbulence community on spontaneous stochasticity and mixing. Let us consider these aspects in turn.

    Positivity of the initial vorticity for the Navier-Stokes equations, and the fact that our bounds are independent of the Reynolds number, is intended to allow a comparison with two fundamental rigorous results on existence of weak solutions to the 2D Euler equations. The existence of weak solutions for positive measure-valued vorticity datum was addressed by Delort in~\cite{Delort1991}. On the other hand, Szekelyhidi showed in~\cite{Szekelyhidi2011} that there exist infinitely many weak solutions starting from the unperturbed interface and that infinitely many such solutions preserve the total kinetic energy.

We are also inspired by the notion of universal properties for turbulent mixing in the Kelvin-Helmholtz instability as recently discussed by Thalabard, Bec and Mailybaev~\cite{Thalabard2020}. They observed that the size of the vortex mixing layer and the energy correlation follow empirical scaling laws. In particular, the mixing layer is observed to grow linearly in time and the authors in~\cite{Thalabard2020} conjectured that this linear growth is universal in the sense that it holds asymptotically for large time for generic initial data.

Our main result,  stated precisely in Theorem~\ref{thm:energy-bounds} below, establishes one side of this conjecture while also being consistent with the rigorous analysis of Delort and Szekelyhidi. We prove that the growth of the mixing layer is at most $O(t)$.  The constant in this bound is explicitly determined and depends solely on the strength of the vortex sheet. Theorem~\ref{thm:energy-bounds} does not include a matching lower-bound. Indeed, such a bound is not possible without making further assumptions on the initial data. The solution to the 2D Navier-Stokes equation from the unperturbed sharp interface is a shear flow with a mixing layer that grows as $O(\sqrt{t})$.

Our work is consistent with~\cite{Delort1991,Szekelyhidi2011} in the following sense: any weak solution to the 2D Euler equations (with the unperturbed interface as initial data) that is obtained as a vanishing viscosity limit of solutions to the Navier-Stokes equations  must also satisfy the bounds on the mixing layer we establish. Further, our analysis relies on a comparison with a coarse-grained profile given by the entropy solution to a scalar conservation law. Loosely speaking, the entropy solution to this scalar conservation law corresponds to averaging (in the sense of weak convergence) over the weak solutions to the Euler equations that conserve kinetic energy constructed in~\cite{Szekelyhidi2011}.

%Our main result is a linear (in $t$) upper bound on the growth of the vortex mixing layer over time in line with the conjecture in~\cite{Thalabard2020} and the results on the Euler equations by Delort~\cite{Delort1991} and Szekelyhidi~\cite{Szekelyhidi2011}. 
Our methods and model also differ from~\cite{Delort1991,Szekelyhidi2011} in several ways. We work with the Navier-Stokes equations, not the Euler equations. The main techniques are adapted from the study of mixing zones for the Saffman-Taylor and Rayleigh-Taylor instability~\cite{kalinin2023scale,OttoMenon2004}. These  papers adapt the background field method of Constantin and Doering~\cite{Doering} to the setting of scale invariant estimates for turbulent mixing. A key idea is to construct a background field as the entropy solution to a scalar conservation law. The use of background fields may itself be dated to the work of Hopf~\cite{Hopf}.

An important difference between the analysis of the Kelvin-Helmholtz instability and the Rayleigh-Taylor and Saffman-Taylor instabilities is the absence of a natural density. This is why we require positivity of the vorticity when adapting the techniques of~\cite{kalinin2023scale,OttoMenon2004} to the Kelvin-Helmholtz instability. This assumption is necessary to have a well-defined length scale for the mixing layer. It is somewhat surprising that this (technical) condition agrees with the assumption required by Delort for the Euler equations~\cite{Delort1991}.

Let us now state the model and our results more precisely.
%Finally, our analysis relies on related to a 

    %but our model and techniques are different. 

%We choose data with positive initial vorticity, so that our bounds may be contrasted with rigorous results on the 2D Euler equations with discontinuous shear flow data. 

%    The bounds and initial data are chosen so that we may contrast mixing in the Navier-Stokes
    
%    Consider a discontinuous shear flow in an inviscid, incompressible fluid with opposing flow directions. This flow profile is known to be unstable: an arbitrarily small perturbation leads to the development of a turbulent mixing layer.  The evolution of such shear flows may be described by weak solutions to the Euler equation with vortex sheet initial datum. 

    %Despite the chaotic nature of the flow, certain characteristics are observed to be universal. 

   % In this paper, we study the solutions selected by the vanishing viscosity limit of the Navier-Stokes equations. We build on ideas from~\cite{OttoMenon2004} and~\cite{kalinin2023scale} to establish uniform upper bounds on the size of the mixing zone in the vanishing viscosity limit. 

    \subsection{The model}
    \label{subsec:equations}
    We consider the two-dimensional Navier-Stokes equations for a velocity field $u = (u^y, u^z)$
    \begin{equation}\label{eq:nse}
    \begin{split}
        \partial_t u + (u \cdot \nabla)u + \nabla p - \Delta u =& \,0, \\ 
        \divr u =& \,0, 
    \end{split}
    \end{equation}   
    in the periodic channel $Q_L = L \mathbb{T} \times \mathbb{R} \ni x = (y, z)$. The vorticity $\omega := \nabla^\perp \cdot u = \partial_y u^z - \partial_z u^y$ satisfies
    \begin{equation}    \label{eq:nse-vorticity-form}
        \partial_t \omega + u \cdot \nabla \omega - \Delta \omega = 0.
    \end{equation}
    We define the following reference shear profile $u_0 e_y$, where $u_0 : \R \rightarrow \R $ with a constant~$U > 0$
    \begin{equation}\label{eq:sing-shear-flow}
        u_0(z) = \begin{dcases}
            -\frac{U}{2}, &z \geq 0, \\
            \frac{U}{2}, & z< 0. 
        \end{dcases}
    \end{equation}
   We also introduce the following notations for an arbitrary function $f: Q_L \to \mathbb{R}$,
    \begin{align}
       &\text{horizontal average} & \Bar{f}(z) & = \frac{1}{L} \int_{L \mathbb{T}} f(y, z) \; dy, \\ 
       &\text{normalized integral}&\dashint f(x) \; dx & = \int_{\mathbb{R}} \Bar{f}(z) \; dz.
    \end{align}
    The initial datum $u(\cdot, 0)$ satisfies
    \begin{equation}\label{eq:initial-datum}
        u(x, 0) = u_0(z)e_y \quad \text{ in }Q_L \cap \{ |z| \geq 1 \},
    \end{equation}
    \begin{equation}    \label{eq:initial-datum-conditions}
        \omega (x, 0) \geq 0, \quad
        \dashint z \omega(x, 0) \; dx = 0, \quad 
        \dashint \omega(x, 0) \; dx = U. 
    \end{equation}
    \begin{figure}
        \centering
        \begin{tikzpicture}[scale = 0.85]

\tikzstyle{every path}=[line width=1pt]
%left frame
\draw  plot[smooth, tension=.7] coordinates {(-2,0) (-7,0)};
\draw  plot[smooth, tension=.7] coordinates {(-7,0) (-7,4.5)};
\draw  plot[smooth, tension=.7] coordinates {(-7,4.5) (-2,4.5)};
\draw  plot[smooth, tension=.7] coordinates {(-2,4.5) (-2,0)};
\draw  plot[smooth, tension=.7] coordinates {(-7,0) (-7,-4.5)};
\draw  plot[smooth, tension=.7] coordinates {(-7,-4.5) (-2,-4.5)};
\draw  plot[smooth, tension=.7] coordinates {(-2,0) (-2,-4.5)};

%right frame
\draw  plot[smooth, tension=.7] coordinates {(2,0) (2,4.5)};
\draw  plot[smooth, tension=.7] coordinates {(2,4.5) (7,4.5)};
\draw  plot[smooth, tension=.7] coordinates {(7,4.5) (7,0)};
\draw  plot[smooth, tension=.7] coordinates {(2,0) (2,-4.5)};
\draw  plot[smooth, tension=.7] coordinates {(2,-4.5) (7,-4.5)};
\draw  plot[smooth, tension=.7] coordinates {(7,0) (7,-4.5)};

%the instability

\draw  plot[smooth, tension=.7] coordinates {(7,0) (6,0) (5,0.5) (4.5,0.5) (4.5,0) (5,-0.5) (5.5,-0.5) (5.5,0)};
\draw  plot[smooth, tension=.7] coordinates {(2,0) (2.5,0) (3.5,-0.5) (4,-0.5) (4,0) (3.5,0.5) (3,0.5) (3,0)};
\draw  plot[smooth, tension=.7] coordinates {(2,-0.5) (2.5,-0.5) (4,-1) (4.5,-1) (4.5,-0.5) (4,0.5) (4,1) (5,1) (6,0.5) (7,0.5)};
\draw  plot[smooth, tension=.7] coordinates {(2,1) (3.5,1) (4,1.5) (6,1) (7,1)};
\draw  plot[smooth, tension=.7] coordinates {(2,-1) (3,-1) (4.5,-1.5) (5,-1) (7,-1)};

%initial data
\draw[<-]  plot[smooth, tension=.7] coordinates {(-6,1.5) (-3,1.5)};
\draw[<-]  plot[smooth, tension=.7] coordinates {(-6,3) (-3,3)};
\node at (-4.5, 2.25) {\large $-U/2$};
\draw[->]  plot[smooth, tension=.7] coordinates {(-6,-1.5) (-3,-1.5)};
\draw[->]  plot[smooth, tension=.7] coordinates {(-6,-3) (-3,-3)};
\node at (-4.5, -2.25) {\large $U/2$};

%arrows

\draw  plot[smooth, tension=.7] coordinates {(4,1.5) (8.5,1.5)};
\draw  plot[smooth, tension=.7] coordinates {(4.5,-1.5) (8.5,-1.5)};
\draw[<->]  plot[smooth, tension=.7] coordinates {(8,1.5) (8,-1.5)};
\node at (8.7 5, 0) {\large $l(t)$};
\draw[->]  plot[smooth, tension=.7] coordinates {(-7.5,0) (-7.5,2.5)};
\node at (-1,-0.5) {\large $y = L$};
\node at (-8,-0.5) {\large $y = 0$};
\node at (-8, 2.5) {\large $z$};
\draw[<-]  plot[smooth, tension=.7] coordinates {(3,3) (6,3)};
\draw[->]  plot[smooth, tension=.7] coordinates {(3,-3) (6,-3)};

\node at (4.5, 2.25) {\large $-U/2$};
\node at (4.5, -2.25) {\large $U/2$};

\end{tikzpicture}
        \caption{Caricature of the instability: The unperturbed shear flow is shown on the left and a crude schematic of the mixing zone of size $l(t)$ is shown on the right. The structure within the mixing zone can be very complicated (see the images in~\cite{VanDyke}). }
        \label{fig:kh-instability}
    \end{figure}

Since $\omega(\cdot,0) \geq 0$ from \eqref{eq:initial-datum-conditions}, we deduce that $\omega(\cdot, t) \geq 0$ for any $t>0$ from the maximum principle. Following~\cite{Thalabard2020}, we consider the measure $d\mu = \Bar{\omega} dz / U \in \mathcal{P} (\mathbb{R})$ and the following two quantities
    \begin{align}
        &\text{expectation} & m(t) =& \int_{\mathbb{R}} z d\mu = \frac{1}{U}\dashint z \omega \; dx, \label{eq:expectation-definition} \\
        &\text{variance} & l^2(t) =& \int_{\mathbb{R}} z^2 d\mu = \frac{1}{U} \dashint z^2 \omega \; dx.\label{eq:variance-definition}
    \end{align}
    
    These are the most natural empirical statistics for the vertical length scale. The center of the mixing layer is given by $m(t)$ and $l(t)$ measures its spreading. As we show later, $m(t)$ is preserved over the time, and taking into account~\eqref{eq:initial-datum-conditions} yields $m(0) = 0$, which means that $l(t)$ is exactly the variance.

    We also introduce two other quantities
    \begin{align}
        &\text{renormalized energy} & E(t) & = \frac{1}{2U} \dashint \frac{U^2}{4} - |u|^2 \; dx, \label{eq:renormalized-energy-definition} \\ 
        &\text{energy separation} & D(t) & = \frac{1}{2U}\dashint |u - u_0 e_y|^2 \; dx \label{eq:separation-definition}.
    \end{align}
    The renormalized energy~\eqref{eq:renormalized-energy-definition} is a measure of dissipation (as shown below) and the energy separation~\eqref{eq:separation-definition} measures the deviation of the perturbed solution from the unperturbed vortex sheet.

\subsection{A remark on the scaling: the role of $L$ as Reynolds number}
The parameter $L$ is the Reynolds number in the following sense.  We rescale the length, time and velocity as follows
\begin{equation}
    \label{eq:rescale1}
    x' = \frac{x}{L}, \quad t'= \frac{t}{L}, \quad u'(x',t') = u(x,t), \quad p'(x',t') = p(x,t).
\end{equation}
Then the domain is the periodic channel $Q = \mathbb{T} \times \mathbb{R}$ with unit width and the Navier-Stokes equations are modified to
\begin{equation}\label{eq:nse2}
    \begin{split}
        \partial_{t'} u' + (u' \cdot \nabla')u' + \nabla' p' - \frac{1}{L} \Delta' u' =& \,0, \\ 
        \divr' u' =& \,0. 
    \end{split}
    \end{equation}  
Thus, we see that the $L\to \infty$ limit is the vanishing viscosity limit (or the infinite Reynolds number limit).  Equation~\eqref{eq:nse2} is the more familiar scaling for the vanishing viscosity limit. On the other hand, it is also clear that the analysis of equation~\eqref{eq:nse2} is equivalent to~\eqref{eq:nse}. When studying mixing layers it is preferable to use the convention~\eqref{eq:nse} and work on the domain $Q_L$ since it also allows us to study transverse structures (such as `fingers' in ~\cite{OttoMenon2004}) while keeping the parameter dependence transparent.

While the rescaled vorticity
\begin{equation}\label{eq:nse3}
    \omega'(x',t') = L \omega(x,t),
    \end{equation}   
is divergent in $L^\infty$, the rescaling ensures precompactness of $\omega'$ in the space of positive measures $\mathcal{M} (Q)$. Indeed, under this rescaling the vorticity measure for the initial profile~\eqref{eq:initial-datum}, converges to the vorticity for the singular shear flow~\eqref{eq:sing-shear-flow} in the weak-$*$ topology on $\mathcal{M} (Q)$. For these reasons, we will study the Navier-Stokes equations in the form~\eqref{eq:nse}--\eqref{eq:initial-datum-conditions}.

    \section{Statement of the main result}\label{sec:main-result}
    \subsection{Uniform estimates on the main quantities} We define a scalar \textit{coarse-grained profile} $\Bar{u} := \overline{u^y}$ for the velocity field $u$. 
    \begin{Theorem}[Growth of the mixing layer]\label{thm:energy-bounds}
        For any $L>0$ and any smooth initial data $u(\cdot, 0)$ satisfying \eqref{eq:initial-datum} and \eqref{eq:initial-datum-conditions}, there exists a unique classical solution to~\eqref{eq:nse} satisfying the following bounds uniformly in $L$
        \begin{equation}\label{eq:length-energy-bound}
            \limsup_{t \to \infty} \frac{l}{Ut} \leq \frac{1}{2\sqrt{3}}, \quad
            \limsup_{t \to \infty} \frac{E}{U^2 t} \leq \frac{1}{12},
        \end{equation}
        In addition, it holds
        \begin{equation}\label{eq:separation-bound}
            \limsup_{t \to \infty} \frac{D}{U^2 t} \leq \frac{1}{4\sqrt{3}},
        \end{equation}
        and $m(t) = 0$ for any $t>0$.
        
        % \begin{equation}\label{eq:length-energy-bound}
        %     \begin{aligned}
        %         \limsup_{t \to \infty} \frac{l}{Ut} & \leq \frac{1}{2\sqrt{3}}, \\ 
        %         \limsup_{t \to \infty} \frac{E}{U^2 t} & \leq \frac{1}{12}, \\
        %         \limsup_{t \to \infty} \frac{D}{U^2 t} & \leq \frac{1}{4\sqrt{3}}.
        %     \end{aligned}
        % \end{equation}
        % \begin{align}
        %     \limsup_{t \to \infty} \frac{l}{Ut} & \leq \frac{1}{2\sqrt{3}},\label{eq:length-bound} \\ 
        %     \limsup_{t \to \infty} \frac{E}{U^2 t} & \leq \frac{1}{12},\label{eq:energy-bound} \\
        %     \limsup_{t \to \infty} \frac{D}{U^2 t} & \leq \frac{1}{4\sqrt{3}}.\label{eq:separation-bound}
        % \end{align}
        % \textcolor{red}{In addition, it holds $m(t) = m(0) = 0$.}
    \end{Theorem}
    \subsection{Relation to scalar conservation laws} The bounds in~\eqref{eq:length-energy-bound} turn into equalities if the coarse grained profile $\Bar{u}$ satisfies the Riemann problem
    \begin{equation}\label{eq:riemann-problem}
        \partial_t \bar{u} - \frac{1}{2} \partial_z \bar{u}^2 = 0, \quad \bar{u}(z, 0) = u_0.
    \end{equation}
    This is reflected in the crucial \textit{scale-invariant} interpolation inequality in \Cref{lm:menon-otto} applied in \Cref{lm:interpolation}. If the coarse-grained profile is indeed given by \eqref{eq:riemann-problem}, the constants in~\eqref{eq:length-energy-bound} are sharp. However, the sharpness of the constant in the separation bound~\eqref{eq:separation-bound} is less clear.

    The entropy solution to the Riemann problem~\eqref{eq:riemann-problem} is a rarefaction wave given by
    \begin{equation}\label{eq:coarse-grained-rarefaction}
        \bar{u}(z) = \begin{dcases}
            U/2, & z \leq -Ut/2, \\ 
            -z/t, &-Ut/2 \leq z \leq Ut/2, \\
            -U/2, & Ut/2 \leq z.
        \end{dcases}
    \end{equation}
    The optimal profile $\bar{u}$ provides a linear in time upper bound on the growth rate of the turbulent mixing layer conjectured in~\cite{Thalabard2020}. The authors of~\cite{Thalabard2020} found a universal asymptotic law, $l = \alpha Ut$  with a pre-factor $\alpha = 0.029$, which does not depend on regularization. \Cref{thm:energy-bounds} implies a rigorous bound on the pre-factor $\alpha$, $\alpha \leq 1/2\sqrt{3} \approx 0.289$. On the other hand, the coarse-grained evolution suggested by~\eqref{eq:riemann-problem} coincides with the conservative subsolutions to the Euler equations constructed in~\cite{Szekelyhidi2011}.

    Let us compute the main quantities for the coarse-grained profile $\bar{u}$. For $l^2$, we have
    \begin{equation}
        \begin{aligned}
            l^2(t) & = \frac{1}{U}\dashint z^2 \bar{\omega} \; dx = \frac{2}{U} \int z(\bar{u} - u_0) \; dz \\ & = \frac{4}{U} \int_0^{Ut/2}z \left( -\frac{z}{t} + \frac{U}{2} \right) \; dz = 2 \int_0^{Ut/2} z \; dz - \frac{4}{Ut} \int_0^{Ut/2} z^2 \; dz \\ & = \frac{U^2 t^2}{4} - \frac{U^2 t^2}{6} = \frac{U^2 t^2}{12}.
        \end{aligned}
    \end{equation}
    For the renormalized energy $E$, it holds
    \begin{equation}
        \begin{aligned}
            E & = \frac{1}{2U} \dashint \frac{U^2}{4} - |u|^2 \; dx = \frac{1}{2U} \int \frac{U^2}{4} - \bar{u}^2 \; dz \\ & = \frac{1}{2U} \int_{-Ut/2}^{Ut/2} \frac{U^2}{4} - \left( \frac{z}{t} \right)^2 \; dz = \frac{U^2 t}{8} - \frac{U^2 t}{24} = \frac{U^2 t}{12}.
        \end{aligned}
    \end{equation}
    The constants above coincide with the constants in~\eqref{eq:length-energy-bound}. Finally, we have for the separation $D$,
    \begin{equation}
        \begin{aligned}
            D & = \frac{1}{2U} \dashint |u - u_0 e_y|^2 \; dx = \frac{1}{2U} \int (\bar{u} - u_0)^2 \; dz \\ & = \frac{1}{2U} \int_{-Ut/2}^0 \left(- \frac{z}{t} - \frac{U}{2} \right)^2 \; dz + \frac{1}{2U} \int_0^{Ut/2} \left( - \frac{z}{t} + \frac{U}{2} \right)^2 \; dz \\ & = \frac{1}{U} \int_0^{Ut/2} \left(-\frac{z}{t} + \frac{U}{2} \right)^2 \; dz = \frac{U^2 t}{24} < \frac{U^2 t}{4 \sqrt{3}}.
        \end{aligned}
    \end{equation}

    % \begin{Remark}      
    %     We contrast the coarse-grained profile $\Bar{u}$ with the entropy solution $s$ to the following Riemann problem for a scalar conservation law 
    %     \begin{equation}\label{eq:riemann-problem}
    %         \partial_t s - \frac{1}{2}\partial_z s^2 = 0, \; s(z, 0) = u_0(z).
    %     \end{equation}
    %     The entropy solution to the Riemann problem~\eqref{eq:riemann-problem} is a rarefaction wave given by
    %     \begin{equation}
    %         s(z) = \begin{dcases}
    %             U/2, & z \leq -Ut/2, \\ 
    %             -z/t, &-Ut/2 \leq z \leq Ut/2, \\
    %             -U/2, & Ut/2 \leq z.
    %         \end{dcases}
    %     \end{equation}
    %     As we shall see below, a crucial interpolation inequality in \Cref{lm:menon-otto} and \Cref{lm:interpolation} turns into an equality if $\Bar{u} = s$. Thus, the entropy solution $s$ provides an upper bound on the turbulent mixing layer that grows linearly in time. 
    % \end{Remark}
    % \begin{Remark}
    %     The bound on the length scale $l$ corresponds to the linear growth of the vortex mixing layer observed and conjectured in~\cite{Thalabard2020}. The coarse-grained evolution suggested by~\eqref{eq:riemann-problem} coincides with the conservative subsolutions to the Euler equations constructed in~\cite{Szekelyhidi2011}.
    % \end{Remark}

    \section{Proofs}\label{sec:proofs}

    The proof of Theorem~\ref{thm:energy-bounds} consists of three parts: (i) the existence and regularity in \Cref{subsec:well-posedness}; (ii) an interpolation inequality related to scalar conservation laws in \Cref{subsec:interpolation}; and (iii) the energy balance relations of the Navier-Stokes equations in \Cref{subsec:energy-balance}. We conclude the proof of \Cref{thm:energy-bounds} in \Cref{subsec:proof-of-theorems}.

    \subsection{Existence of a classical solution}  \label{subsec:well-posedness}
    In this subsection, we assume without loss of generality that $L = 1$ and denote $Q = \mathbb{T} \times \mathbb{R}$.

    \begin{Theorem} \label{t:existence}
        For any smooth initial data $u(\cdot, 0)$ satisfying \eqref{eq:initial-datum} and any fixed $T>0$, there exists a unique smooth solution $u$ to \eqref{eq:nse} in $Q \times [0, T]$. 
    \end{Theorem}
    \begin{proof}[Sketch of the proof]
        Define $v_b(x) := u(x,0) \in C_c^\infty( Q )$ and $v = u - v_b$. From \eqref{eq:nse}, we deduce that $v$ satisfies $v(\cdot, 0) = 0$ and
        \begin{align*}
            \partial_t v + (v \cdot \nabla) v + (v_b \cdot \nabla) v + (v \cdot \nabla) v_b + \nabla p - \Delta v  =& F, \\
            \divr v =& 0, \\ 
            F :=& - (v_b \cdot \nabla) v_b + \Delta v_b.
        \end{align*}
        Then we use Leray's proof reproduced in~\cite{Lemarie2018} to prove the existence, regularity and uniqueness of $v$. The only difference is that the energy estimate for $v$ is given by
        \begin{equation}
                \frac{1}{2} \frac{d}{dt} \int_Q |v|^2 \; dx \leq \int_Q |v|^2 | \nabla v_b| 
                    - \int_Q |\nabla v|^2 \; dx + \int_Q v \cdot F \; dx.
        \end{equation}
        With the above estimate, the existence of one weak solution follows from Theorem 12.2 in ~\cite{Lemarie2018}. By the maximum principle, for $\omega$ satisfying \eqref{eq:nse-vorticity-form}, we deduce that $\omega$ is bounded. Then by standard parabolic theory, the solution is smooth. The uniqueness of $v$, and subsequently $u$, follows from the weak-strong uniqueness in Theorem 12.3 of \cite{Lemarie2018}.

        \begin{Remark}  \label{r:bound_at_infinity}
            By standard Guassian estimates for parabolic equations in \cite{Aronson2017}, we deduce that $\omega$ has exponential decay at $z = \infty$. From standard estimates on higher order derivatives, the vector field $v$ in the proof of \Cref{t:existence} satisfies that, for fixed $T>0$,
            \begin{align}
                v, \nabla v, \nabla^2 v \in& L^2(Q \times [0,T]) \cap L^\infty(Q \times [0,T]). 
            \end{align}
        \end{Remark}
        \noindent Relevant estimates are also proved in \cite{Gallay2013}.
    \end{proof}

    \subsection{An interpolation inequality}  
    \label{subsec:interpolation}

    We state the following interpolation inequality while omitting its proof, since the proof is analogous to Theorem 4 in \cite{OttoMenon2004} or Theorem 3 in \cite{kalinin2023scale}.
    
    \begin{Lemma}[A variation of~{\cite[Theorem~4]{OttoMenon2004}} or ~{\cite[Theorem~3]{kalinin2023scale}}] \label{lm:menon-otto}
        Let $g: [-U/2, U/2] \to \mathbb{R}$, $g \in W^{2, 1}$ be a concave function satisfying the growth condition
        \begin{equation}
            g(\xi) \leq C \xi^{\alpha},\quad\text{with}\quad\alpha > \frac{1}{2}.
        \end{equation}
        If $s: \mathbb{R} \to [-U/2, U/2]$ is a measurable and non-increasing function, then it holds
        \begin{equation}\label{eq:interpolation}
            \int_{\mathbb{R}} g(s(z)) \; dz \leq C_\# \left(\int_{\mathbb{R}} z(s - s_0) \; dz \right)^{1/2}, \quad\text{with}\quad s_0(z) = \begin{dcases}
                +U/2, & z \leq 0, \\
                -U/2, & z > 0.
            \end{dcases}
        \end{equation}
        The sharp constant is given by
        \begin{equation}
            C_\# = \left(2 \int_{-U/2}^{U/2} (g'(s))^2 \; ds  \right)^{1/2}.
        \end{equation}
        The optimal profile $s_g: \mathbb{R} \to [-U/2, U/2]$ for~\eqref{eq:interpolation} is given by
        \begin{equation}
            g'(s_g(z)) = z.
        \end{equation}
    \end{Lemma}
    \begin{Remark}
    The optimal profile above corresponds to the entropy solutions of the Riemann problem
        \begin{equation}
            \partial_t s + \partial_z (g(s)) = 0,\quad s(z, 0) = s_0.
        \end{equation}
    \end{Remark}
    \subsection{Energy balance} \label{subsec:energy-balance} 
    In this subsection, we prove a-priori estimates for the solutions $u$ to the Navier-Stokes equations and their coarse-grained profile $\bar u$.
    
    \begin{Lemma}   \label{lm:energy-growth}
        Let $u$ be a classical solution to~\eqref{eq:nse} with initial data satisfying~\eqref{eq:initial-datum}. Then
        \begin{equation}
           \frac{dE}{dt} = \frac{1}{U} \dashint |\nabla u|^2 \; dx.
        \end{equation}
    \end{Lemma}
    \begin{proof}
        From the proof of \Cref{t:existence} and \Cref{r:bound_at_infinity} we know that $u \in \dot{H}^1$. Also the fact that $u^z \in L^2$ follows from \eqref{eq:initial-datum}. We then compute 
        \begin{equation}
            \begin{aligned}
                \frac{dE}{dt} & = - \frac{1}{U} \dashint u \cdot \partial_t u \; dx \\ & = \frac{1}{U} \dashint u^y( u^y \partial_y u^y + u^z \partial_z u^y) + u^z(u^y \partial_y u^z + u^z \partial_z u^z) \; dx + \frac{1}{U} \dashint |\nabla u|^2 \; dx \\ & = \frac{1}{2U} \dashint u^z \partial_z (u^y)^2 + u^y \partial_y (u^z)^2 \; dx + \frac{1}{U} \dashint |\nabla u|^2 \; dx \\ & = -\frac{1}{2U} \dashint (u^y)^2 \partial_z u^z + (u^z)^2 \partial_y u^y \; dx + \frac{1}{2U}\overline{(u^z (u^y)^2 + u^y (u^z)^2)} \bigg|_{z = -\infty}^{z = +\infty} + \frac{1}{U} \dashint |\nabla u|^2 \; dx \\
                & = \frac{1}{U}\dashint |\nabla u|^2 \; dx.
            \end{aligned}
        \end{equation}
        Here we use the incompressibility condition. We also use \eqref{eq:initial-datum} and \Cref{r:bound_at_infinity} to verify the integrability of various terms in above computation.
    \end{proof}
    
    \begin{Remark}  \label{r:domain_bar_u}
        Since $\omega \geq 0$, we have a pointwise bound $-\partial_z \Bar{u} \geq 0$. Hence, the coarse-grained profile $\Bar{u}$ is monotone in $z$, and satisfies $-U/2 \leq \Bar{u} \leq U/2$ due to the boundary conditions. Indeed, we can compute
        \begin{equation}
            0 \leq \Bar{\omega} = \overline{\partial_y u^z} - \overline{\partial_z u^y} = -\partial_z \Bar{u}.
        \end{equation}
    \end{Remark}
    
    \begin{Lemma}\label{lm:l-derivative}
        Let $u$ be a classical solution to~\eqref{eq:nse} with smooth initial data satisfying \eqref{eq:initial-datum-conditions}, and $\Bar{u}: \mathbb{R} \times \mathbb{R}_+ \to \R $ be its coarse-grained profile. Then
        \begin{align}\label{eq:l-derivative}
            & \frac{dm}{dt} = 0,\quad\text{ and }\\ 
            & l\frac{dl}{dt} = 1 + \frac{1}{U} \int_{\mathbb{R}} g(\Bar{u}(z)) \; dz - E,
        \end{align}
        where $g(s): [-U/2, U/2] \to \mathbb{R}$ is given by
        \begin{equation}\label{eq:flux-definition}
            g(s) = \frac{1}{2} \left(\frac{U^2}{4} - s^2 \right). 
        \end{equation}
    \end{Lemma}
    \begin{proof}
        Both statements follow from the vorticity equation \eqref{eq:nse-vorticity-form} and integration by parts. For the first statement, we have
        \begin{equation}\label{eq:dm-dt-proof}
            \begin{aligned}
                \frac{dm}{dt} & = \frac{1}{U} \dashint z \partial_t \omega \; dx =  \frac{1}{U} \dashint -z \partial_y (\omega u^y) - z \partial_z (\omega u^z) + z \partial^2_y \omega + z \partial^2_z \omega \; dx \\ & = \frac{1}{U} \dashint \omega u^z \; dx - \frac{1}{U} z \overline{\omega u^z} \bigg|_{z = -\infty}^{z = +\infty} - \frac{1}{U}\dashint \partial_z \omega \; dx+ \frac{1}{U} z \partial_z \overline{\omega} \bigg|_{z = +\infty}^{z = +\infty} \\ & = \frac{1}{U} \dashint \omega u^z \; dx. 
            \end{aligned}
        \end{equation}
        Here we have used periodicity in $y$ together with the Gaussian decay of $\omega$ and the boundedness of $|u|$ as $z \to \pm \infty$ in \Cref{r:bound_at_infinity}. Similarly, we compute
        \begin{equation}\label{eq:dm-dt-proof_eq2}
            \begin{aligned}
                \dashint (\partial_y u^z - \partial_z u^y) u^z \; dx & = \frac{1}{2} \dashint \partial_y (u^z)^2 \; dx + \dashint u^y \partial_z u^z \; dx - \overline{u^y u^z} \bigg|_{z = -\infty}^{z = +\infty} \\
                & = -\frac{1}{2} \dashint \partial_y (u^y)^2 \; dx = 0,
            \end{aligned}
        \end{equation}
        where we use the incompressibility condition together with decay of $u^z$ as $z \to \pm \infty$ ($u^z \in L^2$). Then the first statement follows from \eqref{eq:dm-dt-proof} and \eqref{eq:dm-dt-proof_eq2}. 

        To prove the second statement, we again use the vorticity equation \eqref{eq:nse-vorticity-form}. Direct computations give
        \begin{equation}    \label{e:l-derivative:9}
            \begin{aligned}
                l \frac{dl}{dt} & = \frac{1}{2U} \dashint z^2 \partial_t \omega \; dx = \frac{1}{2U} \dashint -z^2 \nabla \cdot (\omega u) + z^2 \Delta \omega \; dx \\ 
                & = 1 + \frac{1}{U} \dashint z \omega u^z \; dx - \frac{1}{2U} z^2 \overline{\omega u^z} \bigg|_{z = -\infty}^{z = +\infty} = 1 + \frac{1}{U} \dashint z\omega u^z \; dx.
            \end{aligned}
        \end{equation}
        Using the definition of $\omega$ yields
        \begin{equation} \label{e:l-derivative:11}
            \dashint z \omega u^z \; dx = \dashint z (\partial_y u^z - \partial_z u^y) u^z \; dx = \dashint z \partial_y \frac{(u^z)^2}{2} - z u^z \partial_z u^y \; dx.
        \end{equation}
        The first term vanishes since $u^z \in L^2$ and $\partial_z u^y \in L^2$ from \Cref{r:bound_at_infinity}. For the second term, integrating by parts yields
        \begin{equation} \label{e:l-derivative:12}
            \begin{aligned}
                \dashint z u^z \partial_z u^y \; dx & = - \dashint u^y \partial_z (z u^z) \; dx + z \overline{u^y u^z} \bigg|_{z = -\infty}^{z = +\infty} \\ & = - \dashint u^y (u^z + z \partial_z u^z) \; dx.
            \end{aligned}
        \end{equation}
        The second term in the second line of \eqref{e:l-derivative:12} vanishes due the incompressibility of $u$ and periodicity in $y$. Then it follows from \eqref{e:l-derivative:9}, \eqref{e:l-derivative:11} and \eqref{e:l-derivative:12} that
        \begin{align}\label{eq:dl-dt-uy-uz}
            \quad l\frac{dl}{dt} = 1 + \frac{1}{U}\dashint u^y u^z \; dx.
        \end{align}

        Now we choose $\Bar{u}e_y$ to be a background field. It is easy to verify $u - \Bar{u}e_y \in L^2(Q)$ and $\overline{\Bar{u} u^z} = 0$ for any $z \in \mathbb{R}$ and $t \geq 0$. Then it holds by~\eqref{eq:dl-dt-uy-uz} and Young's inequality,
        \begin{equation}
            \begin{aligned}
                l \frac{dl}{dt} & = 1 + \frac{1}{U} \dashint u^z(u^y - \bar{u}) \; dx \leq 1 + \frac{1}{2U}\dashint (u^z)^2 + (u^y - \Bar{u})^2 \; dx \\ & \leq 1 + \frac{1}{2U} \dashint (u^z)^2 + (u^y)^2 - 2u^y \bar{u} + (\bar{u})^2 = 1 + \frac{1}{2U} \dashint |u|^2 - (\Bar{u})^2 \; dx \\ & = 1 + \frac{1}{2U} \int_{\mathbb{R}} \frac{U^2}{4} - (\bar{u})^2 \; dz - \frac{1}{2U} \dashint \frac{U^2}{4} - |u|^2 \; dx = 1 + \frac{1}{U} \int_{\mathbb{R}} g(\bar{u}(z)) \; dz - E,
            \end{aligned}
        \end{equation}
        which yields the desired statement.
        
        % Now we choose $\Bar{u}e_y$ to be a background field. It is easy to verify $u - \Bar{u}e_y \in L^2(Q)$ and $\overline{\Bar{u} u^z} = 0$ for any $z \in \mathbb{R}$ and $t \geq 0$. Then it holds 
        % \begin{equation}
        %     \begin{aligned}
        %         l \frac{dl}{dt} & = 1 + \frac{1}{2U}\dashint (u^z)^2 + (u^y - \Bar{u})^2 \; dx \\ & = 1 + \frac{1}{2U}\int_{\mathbb{R}} \frac{U^2}{4} - (\Bar{u})^2 \; dx - \frac{1}{2U}\dashint \frac{U^2}{4} - |u|^2 \; dx \\ & = 1 + \frac{1}{U}\int_{\mathbb{R}} g(\Bar{u}) \; dz - E.
        %     \end{aligned}
        % \end{equation}
    \end{proof}
   
    \begin{Lemma}\label{lm:interpolation}
        Let $u$ be a classical solution to~\eqref{eq:nse} with initial data satisfying~\eqref{eq:initial-datum}, and $\omega$ be the corresponding vorticity. Let $\Bar{u}: \mathbb{R} \times \mathbb{R}_+ \to \R $ be the coarse-grained profile, and $g(\Bar{u})$ be defined by~\eqref{eq:flux-definition}. Then the following inequality holds
        \begin{equation}
            \frac{1}{U} \int_{\mathbb{R}} g(\Bar{u}(z)) \; dz \leq \frac{U l}{\sqrt{12}}.
        \end{equation}
    \end{Lemma}
    \begin{proof}
        The proof is based on applying \Cref{lm:menon-otto}. First, note that it holds
        \begin{equation}    \label{e:interpolation:3}
            \frac{U l^2}{2}  = \frac{1}{2} \dashint z^2 \omega \; dx = \int_{\mathbb{R}} z (\Bar{u} - u_0) \; dz. 
        \end{equation}
        Indeed, by incompressibility, we have
        \begin{equation}    \label{e:interpolation:4}
            \frac{1}{2} \dashint z^2 \omega \; dx = \frac{1}{2} \dashint z^2 (\partial_y u^z - \partial_z u^y) \; dx = -\frac{1}{2} \int_{\mathbb{R}} z^2 \partial_z \Bar{u} \; dz.
        \end{equation}
        Fix $\varepsilon > 0$ and let $u^\varepsilon_0 = u_0 * \theta_\varepsilon$ be a standard mollification. We observe that 
        \begin{equation}    
            \int_{\mathbb{R}} z^2 \partial_z u_0^\varepsilon \; dz \to 0 \text{ as }\varepsilon \to 0,
        \end{equation}
        since $\partial_z u_0 = U \delta_0(z) \in \mathcal{S}'(\mathbb{R})$. Hence we have
        \begin{equation}    \label{e:interpolation:6}
            \begin{aligned}
                & -\frac{1}{2} \int_{\mathbb{R}} z^2 \partial_z (\Bar{u} - u_0^\varepsilon) \; dz - \dashint z^2 \partial_z u_0^\varepsilon \; dx \\ & = \int_{\mathbb{R}} z(\Bar{u} - u^\varepsilon_0) \; dx - \int_{\mathbb{R}} z^2 \partial_z u_0^\varepsilon \; dz \to \dashint z(\Bar{u} - u_0) \; dz. 
            \end{aligned}
        \end{equation}
        Then, combining \eqref{e:interpolation:4} and \eqref{e:interpolation:6} gives \eqref{e:interpolation:3}.

        Finally, \Cref{lm:menon-otto} implies 
        \begin{equation}
            \frac{1}{U} \int_{\mathbb{R}} g(\Bar{u}) \; dz \leq \frac{1}{U} \left(2 \int_{-U/2}^{U/2} s^2 \; ds \right)^{1/2} \frac{U^{1/2}}{\sqrt{2}} l = \frac{U l}{\sqrt{12}},
        \end{equation}
        which concludes the proof.

    \end{proof}

    \subsection{Proof of Theorem~\ref{thm:energy-bounds}}\label{subsec:proof-of-theorems}
    The existence, regularity and uniqueness of the solution are proved in \Cref{t:existence}.
    
    To prove the first bound in Theorem~\ref{thm:energy-bounds}, we combine \Cref{lm:energy-growth}, \Cref{lm:l-derivative} and \Cref{lm:interpolation} to obtain
    \begin{equation}
        l\frac{dl}{dt} = 1 + \frac{1}{U} \int g(\Bar{u}) \; dz - E \leq 1 + \frac{1}{U} \int_{\mathbb{R}} g(\Bar{u}) \; dz \leq 1 + \frac{U l}{\sqrt{12}}.
    \end{equation}
    Solving the differential inequality together, and taking limsup yields a large-time bound on $l$. The second bound follows from Jensen inequality and \Cref{lm:interpolation},
    \begin{equation}
        E \leq \frac{1}{2U} \int_{\mathbb{R}} \frac{U^2}{4} - (\Bar{u})^2 \; dx \leq \frac{U}{\sqrt{12}} \frac{Ut}{\sqrt{12}} = \frac{U^2 t}{12}.
    \end{equation}

    For the last bound in Theorem~\ref{thm:energy-bounds}, we start with the following definition of the separation,
    \begin{equation}
        \begin{aligned}
            \frac{1}{2U} \dashint |u - u_0 e_y|^2 \; dx & = \frac{1}{2U} \dashint \left(|u|^2 - 2\Bar{u}u_0 + |u_0|^2 \right) \; dx = \frac{1}{U} \int_{\mathbb{R}} \left( |u_0|^2 - u_0 \Bar{u} \right) - E.
        \end{aligned}
    \end{equation}

    % \textcolor{red}{This proof relies on the specific form of $\Bar{u}$. I want to avoid it.}
    % \textcolor{gray}{Denote 
    % \begin{equation}
    %     \int_{\mathbb{R}} u_0 (u_0 - \Bar{u}) \; dz = \int_{\mathbb{R}} h(\Bar{u}) \; dz,
    % \end{equation}
    % where $h(\Bar{u})$ is given by
    % \begin{equation}
    %     h(\Bar{u}) = \begin{dcases}
    %         \frac{U}{2} \left(\frac{U}{2} - \Bar{u} \right), & 0 < \Bar{u} \leq \frac{U}{2}, \\
    %         \frac{U}{2} \left( \frac{U}{2} + \Bar{u} \right), & -\frac{U}{2} \leq \Bar{u} \leq 0. 
    %     \end{dcases}
    % \end{equation}
    % Apply \Cref{lm:menon-otto} and the sharp constant $C_\#$ is given by $U$. We obtain
    % \begin{equation}
    %     \frac{1}{U} \int_{\mathbb{R}} u_0 (u_0 - \Bar{u}) \; dz \leq U l \leq \frac{U^2 t}{\sqrt{12}}. 
    % \end{equation}}

    % \textcolor{red}{A robust proof}
    
    \noindent Since $-U/2 \leq \Bar{u} \leq U/2$, we choose $R = l$ and write 
    \begin{equation}
        \begin{aligned}
            \frac{1}{U} \int_{\mathbb{R}} u_0 (u_0 - \Bar{u}) \; dz & = \frac{1}{2} \int_{-\infty}^{-R} \left( \frac{U}{2} - \Bar{u} \right) \; dz + \frac{1}{2} \int_{-R}^0 \left( \frac{U}{2} - \Bar{u} \right) \\ & + \frac{1}{2} \int_0^R \left( \frac{U}{2} + \Bar{u} \right) \; dz + \frac{1}{2} \int_R^{+\infty} \left( \frac{U}{2} + \Bar{u} \right) \; dz \\ & \leq \frac{1}{2R} \int_{\mathbb{R}} z (\Bar{u} - u_0) \; dz + \frac{RU}{2} = \frac{U}{2R} l^2 + \frac{RU}{2} = \frac{Ul}{2}.
        \end{aligned}
    \end{equation}
    The last statement directly follows from \Cref{lm:l-derivative}.

    \section{Acknowledgements}
      The work of GM was supported by NSF grants DMS-2107205 and DMS-2407055. He is also grateful to the the Institute of Mathematics, Academia Sinica, the Max Planck Institute for Mathematics in the Sciences at Leipzig and the University of Oslo for hospitality during the completion of this work. Partial support at Academia Sinica was provided by grant NSTC 113-2115-M-001-009-MY3.

    \bibliographystyle{siam}
    \bibliography{ref}

\end{document}